\documentclass{amsart}

\newtheorem{theorem}{Theorem}[section]

\theoremstyle{definition}
\newtheorem{definition}[theorem]{Definition}

\newtheorem{question}[theorem]{Question}

\newcommand{\noin}{\noindent}
\newcommand{\Lim}{{\rm Lim }\; \varphi }

\theoremstyle{remark}

\numberwithin{equation}{section}

\begin{document}

\title{ Properness of Minimal Surfaces \\ with Bounded Curvature }

\author{G. Pacelli Bessa}
\address{Department of Mathematics, Universidade Federal do Cear\'{a}, Fortaleza, Brazil 60455-760}
\email{bessa@mat.ufc.br}
\thanks{The first author was supported  in part by CAPES-PROCAD GRANT \#0188/00-0.}

\author{Luqu\'{e}sio P. Jorge}
\address{Department of Mathematics, Universidade Federal do Cear\'{a}, Fortaleza, Brazil 60455-760}
\curraddr{}

\email{ljorge@mat.ufc.br}
\thanks{The second author was supported in part by CAPES GRANT \#BEX2067/00-5.}


\date{\today}



\begin{abstract}We show that  immersed minimal surfaces in $\mathbb{R}^{3}$ with bounded curvature and proper self intersections  are  proper. We also show that restricted to  wide components  the immersing map is always proper, regardless the map being proper or not. Prior to these results it was only  known   that injectively immersed minimal surfaces with bounded curvature were proper. 
\end{abstract}

\maketitle

\section{Introduction}

 Most of the results about the structure  of complete minimally immersed    surfaces of $\mathbb{R}^{3}$  requires  the properness of their immersing maps.  It would be interesting question  to ask. What complete geometries  on a  surface  imply that any isometric minimal immersion of the surface into $\mathbb{R}^{3}$  is proper? The first   answering result (to the best of our knowledge) toward this question is due to Rosenberg \cite{rosenberg}. He proved that  a complete    injectively and  minimally immersed  surface of $\mathbb{R}^{3}$ with bounded  curvature is proper\footnote{An immersed surface is said to be proper if its immersing map is proper.} and  presented  an example of a complete  minimal surface dense in all of $\mathbb{R}^{3}$  with bounded curvature. There is also  a beautiful example due to  Andrade \cite{andrade},  of a complete   minimal plane with bounded  curvature, dense in large subsets of $\mathbb{R}^{3}$. These  examples show that the bounds on the curvature alone is not sufficient to ensure properness of an immersed minimal surface.  So, in the class of complete minimal surfaces of $\mathbb{R}^{3}$ with bounded curvature  one has at one hand, that injectively immersed surfaces are proper, and on the other hand there examples of dense (non proper) minimal surfaces.  In the midle of these two extreme cases there are the surfaces that are proper but not injectively immersed. In view of the fact that the failure on properness  under bounded curvature  are the dense  examples, one is tempted to conjecture that unless a complete minimal surface with bounded curvature is dense, it is proper. Although we do not know many examples to support such conjecture, a weaker version of it is true. To state the result precisely  we need the following definition.

\begin{definition} An isometric immersion $\varphi:M\hookrightarrow \mathbb{R}^{3}$ is said to have   proper self intersections    if the restriction of $\varphi $ to $\Gamma = \varphi^{-1}(\Lambda )$ is a proper map, where $\Lambda \, = \,\{x\in \mathbb{R}^{3}; \#\, (\varphi^{-1}(x))\geq2\}$.

\end{definition}

\begin{theorem} \label{thm1} A minimal immersion $\varphi:M\hookrightarrow \mathbb{R}^{3}$ of a complete surface   with bounded curvature and  proper self intersections  is proper.
\end{theorem}
\noin This theorem shows that the key to properness is enconded in the self intersection set $\Gamma $. In some sense it is not necessary too much too have a proper immersion. To complement this result we should look at the conected components of $M\setminus \Gamma$.

\begin{definition} A connected component $M' $ of $M\setminus \Gamma$  is wide if for a divergent sequence of points $x_{k}\in M' $ with $\varphi (x_{k})$ converging in $\mathbb{R}^{3}$  there is a sequence of positive real numbers $r_{k}\rightarrow \infty$ such that the geodesic balls  $B_{M}(x_{k},r_{k})$ of $M$ centered at $x_{k}$ with radius $r_{k}$ are contained in $ M' $. Otherwise we say that the connected component is narrow.
\end{definition}

\begin{theorem}   Let $\varphi:M\hookrightarrow \mathbb{R}^{3}$ be a minimal immersion  of a complete surface with bounded sectional curvature. Then the restriction of $\varphi $ to any wide component $M'$ of $M\setminus \Gamma $ is proper.\label{thm3}
\end{theorem}

 \noin Both results (Theorems \ref{thm1}, \ref{thm3}) are corollaries  of the following theorem.
\begin{theorem}\label{thm2}  Let $\varphi:M\hookrightarrow \mathbb{R}^{3}$ be an isometric minimal immersion  of a complete surface with bounded sectional curvature. There is no divergent sequence of points $x_{k} \in M$ with $dist_{M}(x_{k},\Gamma )\rightarrow \infty$ such that $\varphi (x_{k})$ converges in $\mathbb{R}^{3}$.
\end{theorem}

These results above extend Rosenberg's result in two ways. First, $\Gamma = \emptyset $ if the immersion is injective, thus the restriction of $\varphi$ to $\Gamma $ is trivially proper. Secondly, when the immersion is injective, $M$ is a (the only)  wide component.

The proof of Theorem \ref{thm2} relies the Strong Half-Space Theorem for bounded curvature proved in \cite{bessa-jorge} and in a  well known convergence result. It also can proved using an stability result (Theorem 1.5) also proved in \cite{bessa-jorge}

We shall finish this introduction presenting some questions related to this work that we think it is of importance.

\begin{question}Let $\varphi :M\rightarrow \mathbb{R}^{3}$ be an isometric minimal immersion of a complete surface with bounded curvature. Let $S\subset \Lim $ be a limit leaf. Can $S$ be injectively immersed? Or can $S$ have an injectively immersed end?
\end{question}

 The definitions of $\Lim $ and limit leaf is given below in the preliminaries. The negative answer would be an indicative that $\Lim $ is generated by accumulation points  of $\varphi_{ \mid_{\Gamma}}$. Related to Question 1.6 one can ask,  when $\Lim$ is generated by ${\rm Lim}\varphi_{ \mid_{\Gamma}}$ and when it is not? 

\begin{question} Let $\varphi :M\hookrightarrow \mathbb{R}^{3}$ be a non proper isometric minimal immersion of a complete surface with bounded curvature. Does $M\setminus \Gamma$ have wide connected components?
\end{question}

 The nonexistence of wide components in $M\setminus\Gamma$ would suggest that $M$ is included in the $\Lim .$ Unfortunately there is not enough examples to understand what happens.
Finally, it is important to know when a minimal immersion of $M$ is an O-minimal set or not. O-minimal structures have been studied  in algebraic geometry by many people as L. Br\"oker, M. Coste, L. van den Dries and others, see \cite{coste}. The following question, if answered positively, would imply that a minimal immersion with bounded curvature is an O-minimal set.

\begin{question} Let $\varphi :M\hookrightarrow \mathbb{R}^{3}$ be a non proper isometric minimal immersion of a complete surface with bounded curvature. Does the intersection of $M$ with one line of $\mathbb{R}^3$ have only finite connected components?
\end{question}

 \section{Preliminaries}

 One of the reasons of the  properness requirement in most structure results in minimal surface theory is the absence of tools applicable for non proper minimal immersions. To remedy this situation we introduce the 
  the notion of {\em limit sets} of an isometric immersion. The limit sets of non proper isometric minimal immersions into $\mathbb{R}^{3}$ with bounded curvature  have a rich structure that can be used to better understand those type of immersions.

\begin{definition} Let $\varphi: M\hookrightarrow N$, be an isometric  immersion where $M$ and $N$ are  complete Riemannian manifolds. The set of all points $p\in N$ such that there exists a divergent sequence $\{p_{l}\}\subset M$ so that $\varphi (p_{l})\rightarrow p$ in $N$ is called the limit set  of $\varphi $, denoted by $\Lim$, i.e.
\begin{equation}\Lim \,  = \, \{p\in  N ; \;\exists \,\{p_{l}\}\subset M ,\,{\rm dist }_{M}(p_{0},p_{l})\rightarrow \infty \; {\rm and} \;{\rm dist}_{N}(p, \varphi (p_{l}))\rightarrow 0\}\label{eq1}
\end{equation}
\end{definition}

\noindent The following  convergence theorem is well known in the literature. A more general statement can be proved   for isometric immersions $\varphi :M\hookrightarrow N$ with bounded mean curvature, where $M$ has scalar curvature bounded from below and $N$ has bounded geometry, \cite{bessa-jorge2}.

\begin{theorem}[Folk] Let $\varphi : M\hookrightarrow \mathbb{R}^{3}$ be  a non proper isometric minimal  immersion of a complete surface with bounded sectional curvature. Then  
  for each point $p\in \Lim $ 
there is  a family of complete  immersed minimal surfaces $S_{\lambda}\subset \Lim $, with  bounded sectional curvature, containing $p$, ($p\in S_{\lambda}   $). Each of these surfaces $S_{\lambda}$ is called a limit leaf passing through $p$
\end{theorem}

\begin{proof} The idea of a proof for this theorem is the following: Since the immersing map $\varphi $ is non proper, there is a divergent sequence of points $x_{k}$ in $M$ with $\varphi (x_{k})$ converging to a point $p$ in $\mathbb{R}^{3}$. By the fact that the immersion is minimal, the surface has bounded sectional  curvature and the ambient space bounded geometry, there is a sequence of disjoint disks $D_{k}\subset M$ with radius uniformly bounded from below containing $x_{k}$ such that $\varphi (D_{k})$ are minimal graphs over their tangent spaces at $\varphi (x_{k})$, also with radius uniformly bounded from below. These graphs converge, up to subsequences,  to a minimal graph $D$ with bounded curvature, containing $p$. Each point $q$ of $D$ is limit   of a sequence of points $\varphi (q_{k})\in \varphi (D_{k})$. Repeating this process for the sequence $\varphi (q_{k})$ and $q$ one sees that this limit graph $D$ can be  extended to a complete minimal surface with bounded curvature $S\subset \Lim $ called limit leaf.  One can see that given a compact $K\subset S$, there are a sequence of compacts $K_{k}\subset \varphi (M)$ converging (locally as graphs) to $K$. 
\end{proof}

\section{ PROOFS OF THEOREMS \protect{\ref{thm1}, \ref{thm3}} and \protect{\ref{thm2} } }
\subsection{Theorem \protect{\ref{thm2}}}
\begin{proof}
Suppose we have a divergent sequence of points $x_k \in M$ with $\varphi (x_{k})$ converging to a point $p\in \mathbb{R}^{3}$ and $\mbox{dist}_M(x_k,\Gamma)\rightarrow \infty$. For each $x_{k}$, consider an open disk $D_{k}$ containing $x_{k}$ with radius $\mbox{dist}_M(x_k,\Gamma).$  This sequence converges to a limit leaf    $S$ of $\Lim$ passing by $p$. Suppose (by contradiction) that there is  $q \in S \cap\varphi (M)$, then there is a sequence of points $y_{k}\in D_{k}$ such that $\varphi (y_{k})\rightarrow q$. For large indices, the disks $D(q_{k},1)\subset D_{k}$ centered at $q_{k}$ with radius $1$ are such that $\varphi (D(q_{k},1)) $ are  graphs over a disc in $S$ centered at $q$ with a radius near $1$, thus they intersect $\varphi (M)$ transversally, showing that $D_{k}\cap \Gamma \neq \emptyset $, contradiction. Therefore  $S \cap\varphi (M)=\emptyset$. This leads to another contradiction, because  by  the Strong Half Space Theorem for bounded curvature, $S$ and $\varphi (M)$ are parallel planes. 
\end{proof}

\subsection{Theorems \protect{\ref{thm1} \& \ref{thm3}}}
\begin{proof}
  By Theorem \ref{thm2}, the only choice to  $\varphi $ be non proper is the existence of a divergent sequence $x_{k}\in M$ so that $\varphi (x_{k})$ converges to a point $p\in \mathbb{R}^{3}$ and $\liminf_{k\rightarrow \infty}{\rm dist}(x_{k},\Gamma)<\infty$. But this would be showing that $\varphi_{\mid_{\Gamma}}$ is proper, contradicting the hypothesis of the Theorem \ref{thm1}.
\end{proof}

\begin{proof}

 If $\varphi $ restricted to a wide component $M'\subset M\setminus \Gamma$ is not proper, then there is a divergent sequence $x_{k}\in M'$ such that $\varphi (x_{k})$ converges.  By definition of wide components there is a sequence of positive real numbers $r_{k}\rightarrow \infty$ such that the geodesic balls $B_{M}(x_{k},r_{k})$ of $M$ are contained      in $M'$. Thus  $\liminf_{k\rightarrow \infty}{\rm dist}(x_{k},\partial M')=\infty $.  These geodesic balls $B_{M}(x_{k},r_{k})$ generates a limit leaf $S$ and since the balls are contained in $M'$ the restriction of $\varphi $ to them is injective, thus $S\cap \varphi (M)=\emptyset$. Againd by the Strong Half-Space Theorem for bounded curvature $S$ and $\varphi (M)$ are parallel planes. 

\end{proof}
\bibliographystyle{amsplain}

\end{document}